\documentclass[a4paper]{article}

\usepackage{amssymb}
\usepackage{mathrsfs}

\newtheorem{lemma}{Lemma}
\newtheorem{theorem}{Theorem}
\newtheorem{defin}{Definition}

\newcommand{\scalar}[2]{\langle{#1},{#2}\rangle}
\newcommand{\traccia}{\mathrm{Tr} \:}
\newcommand{\proof}{\textbf{\textit{Proof }}}

\title{An estimate for the entropy of Hamiltonian flows}
\author{F. C. Chittaro}
\date{}

\begin{document}

\maketitle

\begin{abstract}
In the paper we present a generalization to Hamiltonian flows on 
symplectic manifolds of the 
estimate proved by Ballmann and Wojtkovski in \cite{BaWoEnGeo} for the 
dynamical 
entropy of the geodesic flow on a compact Riemannian manifold of nonpositive
sectional curvature. Given such a Riemannian manifold $M,$ Ballmann and 
Wojtkovski proved that the dynamical entropy $h_{\mu}$ of the geodesic flow 
on $M$ satisfies the following inequality:
$$
h_{\mu} \geq \int_{SM} \traccia \sqrt{-K(v)} \; d\mu(v),
$$ 
\noindent
where $v$ is a unit vector in $T_pM$, if $p$ is a point in $M$, $SM$ is
the unit tangent bundle on $M,$ $K(v)$ is defined as $K(v) = \mathcal{R}(
\cdot,v)v$, 
with $\mathcal{R}$ Riemannian curvature of $M$, and $\mu$ is the normalized 
Liouville 
measure on $SM$.

We consider a symplectic manifold $M$ of dimension $2n$,
and a compact 
submanifold $N$ of $M,$ given by the regular level set of a 
Hamiltonian function on $M$; moreover we consider a smooth 
Lagrangian distribution of rank $n-1$ on $N,$ and we assume that
the reduced curvature $\hat{R}_z^h$ of the Hamiltonian vector field $\vec h$ 
is nonpositive. Then we prove that under these assumptions the dynamical 
entropy $h_{\mu}$ of the Hamiltonian flow w.r.t. the normalized Liouville 
measure on $N$ satisfies:
\begin{equation} \label{eq: stima}
h_{\mu} \geq \int_N \traccia \sqrt{-\hat{R}_z^h} \; d\mu.
\end{equation}
\end{abstract}

\section{The curvature}

Let $M$ be a $2n$ dimensional smooth manifold endowed with the symplectic
structure $\sigma.$ Let $h : M \rightarrow  \mathbb{R}$ be a smooth function on the 
manifold, let $\vec h$ denote the Hamiltonian vector field associated 
to it, $d_z h = \sigma (\cdot,\vec h(z))$, and assume that $\vec h$ is a complete
vector field; we will denote by $\phi^t(\cdot) := e^{t\vec h(\cdot)} (\cdot)$ 
the flow generated by $\vec h.$ Let $\Lambda$ be a Lagrangian distribution on $M,$ 
and let us define, for any $z \in M$, the bilinear mapping $g^h_z: \Lambda_z 
\times \Lambda_z \rightarrow \mathbb{R}$ as $g^h_z(X,Y) =
\sigma ([\vec h,X],Y),\; X,Y \in \Lambda_z.$

\begin{defin} 
The Hamiltonian vector field $\vec h$ is said to be regular at $z 
\in M$ w.r.t. the Lagrange distribution $\Lambda$ if the bilinear 
form $g_z^h$ is nondegenerate. A regular Hamiltonian vector field
$\vec h$ is  said to be monotone at $z \in M$ w.r.t. $\Lambda$ if
the form $g_z^h$ is sign-definite.
\end{defin}

\vspace{0.3cm}
\noindent
\textbf{Example} Assume that $\Lambda$ is an involutive Lagrangian
distribution; then, by Darboux-Weinstein Theorem, there exist local 
coordinates $\{(p,q) : p,q \in \mathbb{R}^n \}$ such that 
$\sigma = \sum_{i=1}^n dp_i \wedge dq^i$ and $\Lambda_z  = \{(p,0)\}$;
in these coordinates, the previous requirement about the bilinear form
$g_z^h$ is equivalent to asking the matrix $\{\frac{\partial^2 h}{
\partial p_i \partial p_j}\}$ to be nondegenerate and sign-definite.
\vspace{0.3 cm}

Let us assume that $\vec h$ is regular and monotone.
We define a curve in the Lagrange Grassmannian $L(T_zM)$
putting $J_z(0)= \Lambda_z,\;$ $J_z(t) = \phi^{-t} _* \Lambda_{\phi^t z};$
this curve is called \emph{Jacobi curve}.  
Using the terminology of \cite{AgGaFee}, the curve is regular, 
because the bilinear form $g^h_z$ is nondegenerate;
we have that, for any $t$ sufficiently close (but not equal to) 0, 
$J_z(t)$ is transversal to $J_z(0)$ \cite{AgGeHam}.
Let us denote by $\pi_{J_z(t) J_z(0)}$ the projector of $T_z M$ onto 
$J_z(0)$ and parallel to $J_z(t),$ and note that the space 
$\{\pi_{\Delta J_z(0)} : \Delta \in G_n(T_z M), 
\Delta \in J_z(0)^{\pitchfork}\}$ is an affine subspace of 
$gl(T_z M)$ \cite{AgGeHam}; if we compute the Laurent expansion 
around 0 of the operator-valued function $t \mapsto \pi_{J_z(t) J_z(0)}$,
that is 
$\pi_{J_z(t) J_z(0)} =  \pi_0 + \sum_{i \neq 0} t^i \pi_i,$ we can  
prove that, for $i \neq 0$, $\pi_i \in gl(T_z M)$, while $\pi_0$ is an 
element of the affine space and hence there exists a unique $\Delta \in 
J_z(0)^{\pitchfork}$ such that $\pi_0  = \pi_{\Delta J_z(0)};$ this 
subspace is called the derivative element to $J_z(0)$ and is denoted by 
$J_z^{\circ}(0)$. Analogously, we can apply the same procedure to construct
the derivative element to $J_z(t)$ for $t \neq 0$, and hence we can define
the \emph{derivative curve} of the curve $J_z(t)$: $t \mapsto J_z^{\circ}(t)$;
moreover, we have that $J_z^{\circ}(t) = \phi^{-t}_* J_{\phi^t z}^{\circ} (0).$

Since the Jacobi curve is regular, its derivative curve is smooth and lies 
in the Lagrange Grassmannian of $T_z M$ \cite{AgGeHam}.
These two curves form a splitting (which is called canonical splitting) of 
$T_z M$ into two Lagrangian subspaces $T_z M = J_z(t) 
\oplus J_z^{\circ}(t).$

Let $\Delta_0$ and $\Delta_1$ be two transversal subspaces in the 
Grassmannian $G_n (T_z M),$ and $\xi_0$ and $\xi_1$ be two tangent vectors 
to $G_n(T_z M)$ respectively at the points $\Delta_0$ and $\Delta_1;$ 
let $\gamma_i(t),$ for $i=0,1,$ be two curves in $G_n(T_z M)$ such that
$\gamma_i(0) = \Delta_i$ and $\frac{d}{dt} \gamma_i(t) |_{t=0} = \xi_i.$
Let us set the operator in $gl(\Delta_1):$  
$$[\xi_0,\xi_1] := \frac{\partial^2}{\partial t \partial \tau} \pi_{\gamma_0 (t) 
\gamma_1(0)} \pi_{\gamma_0(0)\gamma_1(\tau)} |_{\Delta_1} \vert_{t=\tau=0};$$

\noindent
this operator depends only on $\xi_0$ and $\xi_1$.

\begin{defin}
The operator $R_{J_z} (t) \in gl (J_z(t))$ defined as 
$$
R_{J_z}(t) := [\dot{J}_z^{\circ}(t),\dot{J_z}(t)]
$$  

\noindent
is called the (generalized) curvature of the curve $J_z(t)$ at the time $t.$ 
\end{defin}

If we choose local coordinates on the Jacobi curve and its derivative curve
putting $J_z(t) \simeq \{(x,S_t x): x \in \mathbb{R}^n\}$ and
$J_z^{\circ}(t) \simeq \{(x,S^{\circ}_t x): x \in \mathbb{R}^n\},$ 
where $S_t$ and $S^{\circ}_t$ are matrices of dimension $n,$ 
the curvature is then $R_{J_z}(t) = (S^{\circ}_t - S_t)^{-1} \dot{
S^{\circ}_t} (S^{\circ}_t - S_t)^{-1} \dot{S_t}.$ 

\begin{defin}
The operator $R_z^h \in gl(J_z(0))$ defined as
$$
R_z^h := R_{J_z}(0)
$$ 
\noindent
is called the curvature of the Hamiltonian vector field $\vec h$ at the point
$z \in M$.
\end{defin}

Let us call $\Sigma_z = \ker (d_z h) / \mathrm{span}
\{\vec h (z)\},$ and let $\psi_z : T_z M \rightarrow T_z M/ \mathrm{span}
\{\vec h(z)\} $ be the canonical projection onto the factor space;
the space $\Sigma_z$ inherits a symplectic structure given by the restriction
of the form $\sigma$.
Let us now set $J_z^h (t) =\phi^{-t}_* [\Lambda_{\phi^t z} \: \cap \: 
\ker (d_{\phi^tz} h) + \mathrm{span}\{\vec h(\phi^t z)\} ]$ (it can be shown that
actually $J_z^h (t) =  J_z (t) \: \cap \: \ker (d_z h) + 
\mathrm{span}\{\vec h(z)\}$), and $\bar{J}_z (t)  = J_z^h (t) / \mathrm{span}
\{\vec h(z)\};$ $\bar{J}_z (t)$ is actually a curve in the Lagrange Grassmannian
$L(\Sigma_z ).$ 
If this Jacobi curve is regular, then its curvature operator $R_{\bar{J}_z} (t)$ is well 
defined on $\bar{J}_z(t).$

\begin{defin}
The operator $\hat{R}_{J_z^h} (t)$ on $J_z^h(t)$ defined as

$$\hat{R}_{J_z^h} (t) := (\psi|_{J_z (t) \cap \ker (d_zh)})^{-1} 
\circ R_{\bar{J}_z} (t) \circ \psi$$
\noindent
is called the curvature operator of the $h-$reduction $J_z^h$ at the time $t$. 
\end{defin}

As before, we define

\begin{defin}
The operator $\hat{R}^h_z$ on $J_z^h(0)$ defined as
$$
\hat{R}_z^h := \hat{R}_{J_z^h} (0)
$$

\noindent
is called the reduced curvature of the Hamiltonian vector field $\vec h$
at the point $z \in M.$
\end{defin}

\vspace{0.3cm}
\noindent
\textbf{Examples} 

\begin{itemize}
\item Let $M = \mathbb{R}^n \times \mathbb{R}^n,$ 
$h(p,q) = \frac{1}{2}|p|^2 + U(q);$ let us consider the Lagrangian distribution
$\Lambda_{(p,q)} = (\mathbb{R}^n,0),$ and let us define 
the Jacobi curve $J_{(p,q)}(t) = \phi^{-t}_* \Lambda_{\phi^t (p,q)};$
Then we have that the curvature is given by $R_{(p,q)}^h =
\frac{\partial^2 U}{\partial q^2},$
and $\hat{R}_{(p,q)}^h  = \frac{\partial^2 U}{\partial q^2} + 
\frac{3}{|p|^2} (\nabla_q U, 0) \otimes (\nabla_q U,0)^T.$
\item Let $M$ be an $n$ dimensional smooth manifold,   
and let $h : T^*M \rightarrow \mathbb{R}$  be such that
the restriction $h|_{T_{\pi(z)}^*M}$ (where $\pi : T^*M \rightarrow M$ is
the canonical projection) is a positive quadratic form, hence it defines a 
Riemannian structure on $M$.
Let $J_z(0) = T_z(T^*_{\pi(z)}M);$ then we have that $R_z^h X = 
\mathcal{R}(\bar{z},\bar{X})\bar{z}$ for any $X \in T_z(T^*_{\pi(z)}M),$
$z \in T^*M,$ where $\mathcal{R}$ is the Riemann curvature tensor, $\bar{z}$ 
is a vector in $TM$ obtained from $z$ by the action of the metric tensor, 
and $X$ is identified with a linear form of $T^*_zM$ via the isomorphism
between $T_z(T^*_{\pi (z)}M)$ and $T^*_{\pi (z)}M.$ The curvature operator of
the $h-$reduction $J_z^h$ is the same, $\hat{R}_z^h = R_z^h$. 
\item Let $M$ as in the previous example, and
let the Hamiltonian function $h$ be the sum of the Hamiltonian function of 
previous example and the function $U \circ \pi,$ where $U$ is a function on
$M$; then
$R_z^h X = \mathcal{R}(\bar{z},\bar{X})\bar{z} + D_{X}(\nabla U)$, and
$\hat{R}_z^h X = R_z^h  X 
+ \frac{3 \scalar{\nabla_{\pi(z)}U}{X}_h}{2 (h(z) - U(\pi(z)))} 
(\nabla_{\pi(z)} U,0)^T,$
where here we denote by $\scalar{\cdot}{\cdot}_h$ the scalar product defined 
by the Riemannian structure given by $h$, and
where $D_X$ is the Riemannian covariant derivative along $X$.
\end{itemize}
\section{Results}

Let $M$ be a $2n$ dimensional smooth manifold endowed with the symplectic 
structure $\sigma$, and let $h: M \rightarrow \mathbb{R}$ a smooth function 
on the manifold; 
we restrict ourselves on a regular sublevel $N$ of the Hamiltonian 
function $h$,
which is then a codimension one submanifold of $M$, and we require this 
submanifold to be compact; moreover, we ask the Hamiltonian function to 
satisfy a regularity condition we will specify later. 
Let us now consider the flow generated by the Hamiltonian vector 
field $\vec h(z),$
and let us notice that it preserves the level sets of
the Hamiltonian, i.e. $h(\phi^t z)=h(z) \; \forall \: t;$
we are interested in computing the dynamical entropy $h_{\mu}(\phi),$ 
where $\mu$ is the (normalized) Liouville measure restricted to the
submanifold $N$; it is defined as $d\mu = \frac{1}{\mathcal{N}} \sigma \wedge
\cdots \wedge \sigma \wedge \iota_X \sigma,$ where $\sigma$ is multiplied by 
itself $n-1$ times, $\iota_X \sigma = \sigma(X,\cdot),$
$X$ is a vector field on a neighborhood of $N$ such that 
$\scalar{dh}{X} =1$ and
$\mathcal{N} = \int_N \sigma \wedge \cdots \wedge \sigma \wedge \iota_X 
\sigma$;
it can be proved that this definition does not depend on the particular 
choice of such a vector field.

In order to compute the dynamical entropy, we are going to use Pesin Theorem
\cite{Mane}, 
which states that the entropy is equal to the integral of the sum of positive 
Lyapunov exponents, 
taken with their multiplicities, and hence we shall compute the exponents
of the Hamiltonian flow. Let us recall that the Lyapunov exponent in the 
point $z \in N$ along the direction $X \in T_z N$ is defined  as 

\begin{equation} \label{eq: lyap}
\lambda^{\pm} (z,X) = \lim_{t \rightarrow \pm 
\infty} \frac{1}{|t|} \log \|\phi^t_* X\|,
\end{equation} 

\noindent
where $\|\cdot\|$ is a scalar product defined on $T_z N$ and, since $N$ is compact,
this definition does not depend on the choice of the norm.  

The symplectic form restricted to $N$ has a one dimensional kernel given by 
the span of the Hamiltonian vector associated to $h:$ indeed 
$$\forall \  v \in T_zN \quad \sigma(v,\vec h) = \scalar{d_zh}{v}=0,$$
since $T_zN=\ker(d_zh);$ 
hence, $\:\forall\: z \in N,$ we can write $T_zN \simeq \Sigma_z \oplus 
\mathrm{span}\{\vec h(z)\},$ where $\Sigma_z = T_z N / \mathrm{span}\{
\vec h(z)\} $ is a $2n-2$ dimensional vector 
space and the restriction $\bar{\sigma} = \sigma |_{\Sigma_z}$ induces 
a symplectic structure on $\Sigma_z$. 
Since $\mathrm{span}\{\vec h\}$ is preserved by the action of its flow, i.e.
$\phi^t_* \vec h(z) = \vec h(\phi^t z),$ we can take the quotient and 
study the exponential divergence of the trajectories along directions given by 
vectors lying in $\Sigma_z$, so we will consider the map 
${\tilde{\phi}^t}_* : \Sigma_z \rightarrow \Sigma_{\phi^t z},$ 
where $ {\tilde{\phi}^t}_*  = {\phi^t}_*|_{\Sigma_z}.$ 

Now we can state the result:


\begin{theorem} \nonumber
Let N be a compact regular level set of a smooth Hamiltonian 
function defined on a smooth symplectic manifold on dimension $2n$; let 
$\Lambda$ be a Lagrangian distribution in $TN/\mathrm{span}\{\vec h\}$ 
and let the Hamiltonian vector field $\vec h$ be monotone on $N$ w.r.t. $\Lambda$.
Consider the Jacobi curve $\bar{J}_z(t)= \tilde{\phi}^{-t}_* \Lambda_{\phi^t z}$ and 
let the curvature $\hat{R}_z^h$ of $\vec h$ be nonpositive.
Then the dynamical entropy $h_{\mu}$ of the Hamiltonian flow on $N$ w.r.t.
the normalized Liouville measure on $N$ satisfies
$$
h_{\mu} \geq \int_N \traccia \sqrt{-\hat{R}_z^h} \;d\mu.
$$ 
\end{theorem}

\noindent
\proof


\noindent
Due to sign-definiteness of the bilinear form $g_z^h$, we can endow 
$\Sigma_z$ with a scalar product; indeed, let us define (for $g_z^h$
positive-definite) the following scalar product on $\bar{J}_z(0):$
$$
\bar{J}_z(0) \ni X,Y \mapsto \scalar{X}{Y}^{\prime}_h:=\bar{\sigma} ([\vec h,X],Y).
$$
\noindent
By means of the symplectic form we can establish an isomorphism between
$\bar{J}_z^{\circ}(0)$ and the dual of $\bar{J}_z(0):$ $\bar{J}_z^{\circ}(0) 
\ni W \mapsto \bar{\sigma}(W,
\cdot) : \bar{J}_z(0) \rightarrow \mathbb{R};$ since there exists a unique $X_W
\in \bar{J}_z(0)$ such that $\bar{\sigma}(W,\cdot)=\scalar{X_W}{\cdot}_h,$ we can define 
the scalar product on $\bar{J}_z^{\circ}(0)$  in this way:
$$
\bar{J}_z^{\circ}(0) \ni W,V \mapsto \scalar{W}{U}^{\circ}_h:=\scalar{X_W}{X_V}_h.
$$

\noindent
Now it is possible to define a scalar product on the whole $\Sigma_z:$ for
any $X,Y \in \Sigma_z$, we set 
$$
\scalar{X}{Y}_h := \scalar{\pi_{\bar{J}^{\circ}_z(0) 
\bar{J}_z(0)}X}{\pi_{\bar{J}^{\circ}_z(0) \bar{J}_z(0)} Y}^{\prime}_h + 
\scalar{\pi_{\bar{J}_z(0) \bar{J}^{\circ}_z(0)} X}{\pi_{\bar{J}_z(0) \bar{J}^{
\circ}_z(0)} Y}^{\circ}_h; $$ 
by definition, $\bar{J}_z^{\circ}(0)$ is orthogonal to $\bar{J}_z(0)$ with respect to the scalar 
product just defined.
%


Since the space $\Sigma_z$ has a symplectic 
structure and for any $t$ the pair $(\bar{J}_z(t),\bar{J}_z^{\circ}(t))$ 
forms a splitting of
Lagrangian subspaces, given a basis $\{\epsilon^1,\ldots,\epsilon^{
n-1}\}$ of $\bar{J}_z(0)$ there is a unique way to choose a basis $\{e_z^1(t),
\ldots,e_z^{n-1}(t)\}$ of $J_z(t)$ such that  $e_z^i(0)=\epsilon^i \; 
\forall\: \: i=1,\ldots,n ,$
$\{\dot{e}_z^1(t),\ldots,\dot{e}_z^{n-1}(t)\}$ is a basis for $J_z^{\circ}(t)$ 
and $\{e_z^i(t),\dot{e}_z^i(t)\}_{i=1}^{n-1}$ is a Darboux basis 
for $\Sigma_z$,
and it is called the \emph{canonical moving frame} \cite{AgGeHam}. 
Moreover, as shown in \cite{AgGeHam}, the vectors $\ddot{e}_z^i(t)$ lie in 
$\bar{J}_z(t)$ for any $i=1,\dots,n-1$, and 
$$
\ddot{e}_z^i(t) = \sum_{j=1}^{n-1} (-R_z(t))_{ij} e_z^j(t),
$$
\noindent
where $R_z(t)$ is the representation of the curvature $\hat{R}_z^h$ w.r.t. 
the basis $\{e_z^i(t)\}_{i=1}^{n-1}$, and it is symmetric.

Let us define, for any $z \in N,$ the basis $\varepsilon_1(z),\ldots,
\varepsilon_{2n-2}(z)$ 
of $\Sigma_z$ by putting $\varepsilon_i(z) = e_z^i(0),
\;\: \varepsilon_{i-n+1}(z) = \dot{e}_z^i(0), \; i=1,\ldots,n-1;$
this basis is indeed orthonormal for any $z.$ 
Consider a vector $X \in  \Sigma_z:$

\begin{equation} \label{eq: base}
 X  =  \sum_{i=1}^{2n-2} x_i \, \varepsilon_i(z) =
  \sum_{i=1}^{n-1}  \eta_i (t) \, e_z^i(t) + \xi_i(t) \, \dot{e}_z^i(t),  
\end{equation}

\noindent
($(\eta_i(t),\xi_i(t))$ are the components of the vector w.r.t. the canonical
moving frame, and obviously $(\eta(0),\xi(0))=(x_1,\ldots,x_{2n-2})$). 
By computations we can prove that the pair $(\eta (t), \xi(t))$
satisfies the differential first-order system

\begin{equation} \label{eq: system}
\Big\lbrace \begin{array}{lcl}  
\dot{\xi}(t)  &=& -\eta(t) \\ 
\dot{\eta}(t) &=& R_z(t)\xi(t)  
\end{array} 
\end{equation}

\noindent
and hence the vector $\xi(t)$ satisfies the second order differential equation

\begin{equation} \label{eq: jacobi}
\ddot{\xi}(t) + R_z(t) \xi(t) = 0.
\end{equation}

\noindent
Since the canonical moving frame is defined such that $e_{\phi^t z}^i(0) =
{\tilde{\phi}^t}_* e_z^i(t),\;$ $\dot{e}_{\phi^t z}^i(0) ={\tilde{\phi}^t}_* 
\dot{e}_z^i(t),\;$ $i=1,\ldots,n-1,$ it implies that $e_z^i(t) = 
{\tilde{\phi}^{-t}} _* e_{\phi^t z} (0) = {\tilde{\phi}^{-t}}_* 
\varepsilon_i(\phi^t z),\;$ $i=1,\ldots,n-1,$ and $\dot{e}_z^i(t) = 
{\tilde{\phi}^{-t}}_* \dot{e}_{\phi^t z} (0) = {\tilde{\phi}^{-t}}_* 
\varepsilon_i(\phi^t z), \;$ $i=n,\ldots,2n-2.$ Hence

\begin{eqnarray}
{\tilde{\phi}^t}_* X & = & \sum_{i=1}^{n-1}  \eta_i (t)\, {\tilde{\phi}^t}_* 
e_z^i(t) + \xi_i(t) \, {\tilde{\phi}^t}_* \dot{e}_z^i(t) \nonumber\\
 & = & \sum_{i=1}^{n-1}  \eta_i \, (t) \varepsilon_i(\phi^t z) + \xi_i(t) \, 
\epsilon_{i+n-1}(\phi^t z) \nonumber,
\end{eqnarray}

\noindent
and it means that the components of ${\tilde{\phi}^t}_* X$ w.r.t. the basis 
$\{\varepsilon_i(\phi^t_z)\}_{i=1}^{2n-2}$ of $\Sigma_{\phi^t z}$ are the same 
as the components of $X$ w.r.t. the canonical moving frame at time $t$. 
Since the basis $\{\varepsilon_i(z)\}_i$ is orthonormal for any $z,$ we find 
that $\|\pi_{\bar{J}_{\phi^tz} \bar{J}^{\circ}_{\phi^tz}} \tilde{\phi}^t_*X\|= 
|\xi(t)|$ and  $\|\pi_{\bar{J}^{\circ}_{\phi^tz} \bar{J}_{\phi^tz}} 
\tilde{\phi}^t_*X\|= |\dot{\xi}(t)|.$


\noindent
Now we shall compute the Lyapunov exponents on $N$; by Multiplicative Ergodic
Theorem \cite{Mane} we know that the limit (\ref{eq: lyap}) exists a.e. (w.r.t. the
standard Liouville measure normalized on $N$) in $N.$ Hence we can define the 
following subspaces of $\Sigma_z:$
\begin{eqnarray}
E^u_z & = & \{X \in \Sigma_z : \lambda^-(z,X) <0\}, \nonumber\\
E^s_z & = & \{X \in \Sigma_z : \lambda^+(z,X) <0\}, \nonumber\\
E^0_z & = & \{X \in \Sigma_z : \lambda^-(z,X) \leq 0 \; \mathrm{and} \; 
\lambda^+(z,X) \leq 0\}; \nonumber
\end{eqnarray}

\noindent
these subspaces span $\Sigma_z$.
For any subspace $E_z$ of $\Sigma_z$ such that $E^u_z \subset E_z \subset 
E^u_z \oplus E^0_z$, we have that $\lim_{t \rightarrow \pm \infty} \frac{1}
{|t|} \log |\det({\tilde{\phi}^t}_*  |_{E_z})|= \pm \chi(z),$
where $\chi(z)$ is the sum of the positive Lyapunov exponents in $z$, taken 
with their multiplicities. 

Knowing this, we are now looking for such a subspace $E_z;$ we'll see that a 
good candidate will be the graph of a proper linear operator, that we will 
call $U_z$, defined from $\bar{J}_z^{\circ}(0)$ to $\bar{J}_z(0)$ .

Let us now introduce for any $z\in N$ the subset $H(z)$ of $\Sigma_z$ 
such that
$$
H(z)=  \{X \in \Sigma_z : \; \frac{d}{dt}\| \pi_{\bar{J}_{\phi^t z}(0) 
\bar{J}^{\circ}_{\phi^t z}(0)} \tilde{\phi}^t_* X\| \geq 0\: \forall\: t\};
$$ 

\noindent
clearly $H(z)$ is intrinsically defined and it is invariant along the 
trajectory $\phi^t z$.

In the following, for simplicity we will denote $\bar{J}_{\phi^tz}(0)$ by $v(t)$ 
and  $\bar{J}^{\circ}_{\phi^tz}(0)$ by $v^{\circ}(t).$

\begin{lemma} \label{lemma: sottospazio}
$H(z)$ is a subspace of $\Sigma_z$.
\end{lemma}

\noindent
\proof From the convexity of $\|\pi_{v(t) v^{\circ}(t)} 
\tilde{\phi}^t_* X\|^2$ (\ref{eq: jacobi}) we deduce that a vector 
$X \in \Sigma_z$ belongs to $H(z)$ if and only if  $\|\pi_{v(t) v^{\circ}(t)}
 \tilde{\phi}^t_* X\|$ is bounded for negative times. Linear combinations
of vectors having this property satisfy this requirement. \hfill $\square$\\

\begin{lemma} $H(z) \cap v(0) = \{0\} $.
\label{lemma: verticale}
\end{lemma}

\noindent
\proof A vector $X \in H(z)$ belongs to $v(0)$ if $\pi_{v(0) v^{\circ}(0)} X=0,$
i.e. if $\xi(0) = 0$; suppose by 
contradiction that such a (nonzero) vector is contained in $H(z);$ then $\frac{
d^2}{dt^2}|\xi(t)|^2|_{t=0}= \scalar{\dot{\xi}(0)}{\dot{\xi}(0)} - 
\scalar{R_z(t)\xi(0)}{\xi(0)} >0,$ hence 0 is a strong minimum for $|
\xi(t)|,$ which contradicts the definition of $H(z).$  \hfill 
$\square$\\

\begin{lemma} $H(z)$ is a Lagrangian subspace. 
\label{lemma: H_Lagrangiano}
\end{lemma}

\noindent
\proof Let us define, for any $\tau \in \mathbb{R}$, $H_{\tau} = \{X \in \Sigma_z :
\frac{d}{dt}\|\pi_{v (t) v^{\circ}(t)} \tilde{\phi}^t_* X \| 
\geq 0 \; \forall \:t \geq \tau\};$  we have that 
$H_{\tau_1} \subseteq H_{\tau_2}$ if $\tau_1 \leq \tau_2$ and that
$H(z)= \cap_{\tau} H_{\tau}.$ 

$H_{\tau}$ contains a Lagrangian subspace for any $\tau.$ 
Indeed, fix $\tau$ and consider $V_{\tau} = \{X \in \Sigma_z : 
\pi_{v (0) v^{\circ} (t)} \tilde{\phi}^{\tau}_* 
X =0 \};$
we prove using coordinates that this subspace is contained in $H_{\tau}:$
if we write $X = \sum_{i=1}^{n-1} -\dot{\xi}(t) e_z^i(t) + \xi(t) 
\dot{e}_x^i(t),$ we have that $\|\pi_{v (t) v^{\circ}(t)} 
\tilde{\phi}^t_* X \| = |\xi(t)|$ and hence, since 
$$
\frac{d}{dt} |\xi(t)|^2|_{t=\tau} =0 \quad \mathrm{and} \quad \frac{d^2}{dt^2}|
\xi(t)|^2 \geq 0 \; \forall\:t,
$$
\noindent
$\frac{d}{dt} \|\pi_{v (t) v^{\circ}(t)} \tilde{\phi}^t_* X \| 
\geq 0 \; \forall\: t \geq \tau,$ and $V_{\tau} \subset H_{\tau}.$

Now, since $\tilde{\phi}^{\tau}_* V_{\tau} = \bar{J}_{\phi^{\tau} z} (0),$ 
and this last subspace is Lagrangian, we proved our claim.
$H(z)$ contains a Lagrangian subspace too;
indeed, let us define for any $\tau$ 
$\hat{H}_{\tau} = \{V \in L(\Sigma_z) :V \subset H_{\tau} \},$ 
which is a compact nonempty subset in the
Lagrange Grassmannian $L (\Sigma_z)$. Moreover, since $\hat{H}_{\tau_1} \subseteq 
\hat{H}_{\tau_2}$ for $\tau_1 \leq \tau_1,$ we have that $\cap_{\tau} \hat{H}_{\tau}
\neq \varnothing;$ hence, since $\hat{H}_{\tau} \subset H_{\tau}$ for any $\tau,$ 
we can conclude that $H(z) \supseteq \cap_{\tau} \hat{H}_{\tau} \neq 
\varnothing,$ that means that $H(z)$ contains a Lagrangian subspace. 

From Lemma \ref{lemma: verticale} we know that $\dim H(z) \leq n-1,$ 
hence we can conclude that $H(z)$ is a Lagrangian subspace. \hfill $\square$\\

Since the space $H(z)$ is Lagrangian and $H(z) \cap \bar{J}_z(0)=0 \; 
\forall\: z,$ there exists a symmetric
linear operator $U_z : \bar{J}_z^{\circ}(0) \rightarrow \bar{J}_z(0) $ 
such that for any element $X \in H(z)$ we have that $X= x+ U_z(0) x,$
where $x \in \bar{J}_z^{\circ}(0)$, i.e. $H(z)$ is the 
graph of the operator $U_z$.    
%
%
%
%
%

Hence we can find a linear operator $V_z : \mathbb{R}^{n-1} \rightarrow
\mathbb{R}^{n-1}$ such that if $H(z) \ni X = \sum_{i=1}^{n-1} \eta_i(0)
\varepsilon_i(z)+\xi(0) \varepsilon_{i+n-1}(z),$ then $\eta(0) 
= -V_z \xi(0),$ and, by (\ref{eq: system}) we get that $\dot{\xi}(0)= 
V_z\xi(0),$ by (\ref{eq: jacobi}) that the operator satisfies 
the equation 
\begin{equation} \label{eq: riccati}
\dot{V}_{\phi^tz} + V_{\phi^tz}^2 + R_z(t)=0.
\end{equation}

\noindent
By definition of $H(z),$ the operator $V_z$ is nonnegative definite for 
any $z$. 

\begin{lemma} $E^u_z \subset H(z) \subset E^u_z \oplus E^0_z.$ 
\label{lemma: contenuto}
\end{lemma}

\noindent
\proof Let $X \in E^u(z)$ and $Y \in E^u_z \oplus E^0_z;$ $\lim_{t 
\rightarrow -\infty} 
\frac{1}{|t|} \log|\bar{\sigma}({\tilde{\phi}^t}_* X,{\tilde{\phi}^t}_*  Y)| 
\leq  \lim_{t \rightarrow -\infty} 
[\frac{1}{|t|} \log\|\bar{\sigma}\| + \frac{1}{|t|} \log\| {\tilde{\phi}^t}_* X\| 
+  \frac{1}{|t|} \log\|{\tilde{\phi}^t}_* Y\|] = \lambda^-(z,X) + 
\lambda^-(z,Y) < 0,$ and this 
implies that $\bar{\sigma}({\tilde{\phi}^t}_* X,{\tilde{\phi}^t}_* Y) 
\rightarrow 0$ for $t \rightarrow 
-\infty .$ Since $\tilde{\phi}^{t*} \bar{\sigma} = \bar{\sigma},$ we get that 
$\bar{\sigma}(X,Y) = 0$ and 
hence $E^u_z$ and $E^u_z \oplus E^0_z$ are skew-orthogonal. By dimensional 
computations, we can prove that actually $E^u_z$ and $E^u_z \oplus E^0_z$ are
the skew-orthogonal complement to each other.

Let $X \in E_z^u$, i.e. $\lim_{t \rightarrow -\infty} \frac{1}{|t|} 
\log\| \tilde{\phi}^t_*X\|< 0,$ and this means that $\| \tilde{\phi}^t_*
X\| <1,$  which implies that $\tilde{\phi}^t_* X$ is bounded in norm for 
nonpositive times, and consequentely also $\pi_{v(t) v^{\circ}(t)} \tilde{\phi}^t_*
X$ is, which implies that $\tilde{\phi}^t_* X \in H(\phi^t z)=
\tilde{\phi}^t_*[H(z)] \Rightarrow E_z^u(z) \subset H(z).$
Since  $H(z)$ is Lagrangian, we also find that  $H(z) \subset 
E_z^u(z) \oplus E_z^0.$  \hfill $\square$\\

\begin{lemma}
Let $X \in H(z);$ then $\pi_{v(0)v^{\circ}(0)} X \in \ker U_z$ if and only 
if $\|\pi_{v^{\circ}(t)v(t)} \tilde{\phi}^t_* X\|=0$ for any $t\leq 0$.
%
%
\end{lemma}

\noindent
\proof We are proving it in coordinates. Let $X$ as in (\ref{eq: base}) 
such that
$\xi(0) \in \ker V_z,$ i.e. $\eta(0)=0;$ since 
by convexity (\ref{eq: jacobi}) $\frac{d^2}{dt^2} |\xi(t)|^2 \geq 0,$ 
and by hypothesis $\frac{d}{dt} |\xi(t)|^2|_{t=0} = 0,$
we get that  
$|\xi(t)|^2$ shall remain constant  $\forall \: t\leq 0,$ which implies, 
using again convexity, that
$|\dot{\xi}(t)|= 0 \ \forall \: t \leq 0.$ 
Conversely, if $\dot{\xi}(t) =0 \ \forall \: 
t\leq 0,$ then obviously we get the thesis.  
\hfill $\square$\\

%
%
Let us denote by $H_0(z)$ the graph of $U_z$ restricted to the 
orthogonal complement in $\bar{J}_z^{\circ}(0)$ to $\ker U_z;$ 
it follows from the above lemma that 
$\tilde{\phi}^t_*[H_0(z)] \subseteq H_0(\phi^t z) \; \forall\:t\geq 0.$ 
Indeed, let $X \in H(z)$ such that 
$\xi(0) \in \ker V_z;$ then, by previous results, $\dot{\xi}(t) =0$ for any
$t \leq 0,$ that means that $\xi(t) \in \ker V_{\phi^tz}$ for negative times,
i.e. $\pi_{v^{\circ}(t) v(t)} \tilde{\phi}^t_* X \in \ker U_{\phi^tz}.$

Since the dimension of $H_0(z)$ is nondecreasing along the orbits of the 
Hamiltonian flow, we get that $\dim H_0(z)$ is constant on a 
$\phi^t-$invariant set of full measure, and hence on this set $\phi^t_*[
H_0(z)]=H_0({\phi^tz}).$

We will work in the space $H_0(z)$ because we need the operator
$U_z$ to be strictly positive definite,
and we are calling $U_z^0$ the restriction of $U_z$ on the orthogonal 
complement to $\ker U_z$ in $\bar{J}_z^{\circ}(0),$ and respectively 
$V^0_z$ and $R_z^0(t)$
the restrictions of $V_z$ and $R_z(t)$ to the orthogonal complement of 
$\ker V_z$ in $\mathbb{R}^{n-1}$; to do this, we shall prove that 
actually it satisfies Lemma \ref{lemma: contenuto}.
First, we need the following result:
\begin{lemma}
$R_z(t)$ vanishes on  $\ker V_{\phi^t z}$ and both $R_z(t)$ and  $V_{
\phi^t z}$ preserve the orthogonal complement in $\mathbb{R}^{n-1}$ to 
$\ker V_{\phi^t z}.$
\end{lemma}

\noindent
\proof 
Call $\Delta_z(t)$ the orthogonal complement in $\mathbb{R}^{n-1}$ to 
$\ker V_{\phi^t z}.$
Let $X \in H(z),$ let $(-\dot{\xi}(t),\xi(t))$ be its 
components as in (\ref{eq: base}), and let $\xi(t) \in \ker V_{\phi^t z}$; 
then, by previous lemma, $\dot{\xi}(\tau) = 0$ for $\tau \leq t,$
which implies the vanishing of the second derivative too, i.e. $R_z(t) 
\xi(t)= 0.$
Let now $x \in \ker V_{\phi^t z}, x' \in \Delta_z(t);$ since $\scalar{x}{
R_z(t) x'} = \scalar{R_z(t)x}{x'}=0,$ we conclude that $R_z(t)[\Delta_z(t)] 
\subseteq \Delta_z(t).$ In the same way we can show that $V_{\phi^t z}
[\Delta_z(t)] \subseteq \Delta_z(t).$ \hfill $\square$\\

Let $X \in H(z) \setminus H_0(z)$; then, $\tilde{\phi}^t_* X$ is 
constant in norm w.r.t. $t$ for any nonpositive $t$.
Hence $\lambda^-(z,X) = 0 \Rightarrow X \notin E^u_z,$ which implies that 
$E^u_z \subset H_0(z).$
Moreover, consider $X = X^{(1)} + X^{(2)} \in H(z);$  we call as usual 
$(-\dot{\xi}^{(i)}(t), \xi^{(i)}(t))$ the components of $\tilde{\phi}^t_*
X^{(i)}$ w.r.t. the orthonormal frame $\{\varepsilon_i(\phi^tz)\}_i$, and we 
assume that
$\xi^{(1)}(t) \in \ker V(z)$ and
$\xi^{(2)}(t)$ lies in $\Delta_z(t)$ (defined as above).
By previous results, we get that $\dot{\xi}^{(1)}(t) = 0 $ for $ t \leq 0,$ 
and hence $\ddot{\xi}^{(1)}(t) = 0 $ for $ t \leq 0,$  and also 
$R_z(t)\xi^{(1)}(t) =0,$ which implies
that both $\xi^{(1)}$ and $\xi^{(2)}$ satisfy equation (\ref{eq: jacobi}).



\noindent
Since $H_0(z)$ is the graph of the operator $U_z^0,$ we 
can express the scalar
product on $H_0(z)$ in term of the scalar product on $\mathbb{R}^{n-1},$
putting $\scalar{X}{Y}_h = \scalar{\xi^X(0)}{A_z(0) \xi^Y(0)}_c,$
where $A_z(t) = \mathbb{I} + {V^0_{\phi^t z}}^2$ ($X$ and $Y$ as above),
and $\scalar{\cdot}{\cdot}_c$ denotes the canonical scalar product on 
$\mathbb{R}^{n-1}.$

We call $a_z(t) = |\det \phi^t_*|_{H_0(z)}|$ the determinant w.r.t. the scalar 
product defined by $A_z(t)$ of $\phi^t_*;$ hence we have that 
$$a_z(t) = \sqrt{ \det A_z(t)} |\det \phi^t_*|_{H_0(z)}|_c = \sqrt{ \det A_z(t)} 
|\det e^{\int_0^t V_{\phi^s z}ds}|_{H_0(z)}|_c .$$

We define $r_z(t) := \frac{d}{dt}\log a_z(t) = \frac{1}{2} \traccia 
\dot{A}_z(t) A^{-1}_z(t) + \traccia V^0_{\phi^t z},$ and we get by 
computations that $r_z(t)= \traccia [(V^0_{\phi^t z} - R_z^0(t) V^0_{\phi^t 
z})(\mathbb{I} + {V_{\phi^t z}^0}^2)^{-1}].$ Since 
$$\chi(z) = \lim_{t \rightarrow \infty} \frac{1}{t} \log |\det(\phi^t_*z|_{
H_0(z)})| = \lim_{t \rightarrow \infty} \frac{1}{t} \log a_z(t) = 
\lim_{t \rightarrow \infty} \frac{1}{t} \int_0^t r_z(s) \: ds,$$

\noindent
by Birkhoff Ergodic Theorem \cite{Mane} we get that, provide that 
$r_z$ is an integrable function on $N,$  
$h_{\mu}(\phi) = \int_N \chi(z) \, d\mu (z) = \int_N r_z (0) \, d\mu.$

Now we are going to compute dynamical entropy using a different
scalar product on $H_0(z),$ after showing that we will get the same value.  
Call $A_z'(t) = V^0_{\phi^t z},$ and define the scalar product 
$\scalar{X}{Y}' = \scalar{\xi^X(0)}{A_z'(0) \xi^Y(0)};$
we also get that $r_z'(t) = \frac{1}{2}\traccia [V^0_{\phi^t z} - R_z^0(t){
V_z^0}^{-1}].$

The volume element on 
$N$ w.r.t. the scalar product given by $A'$ is 
related to the standard volume element in this way: $d\mu' = \sqrt {\frac{\det 
A'}{\det A}} d\mu.$ If we call $c(t) = \frac{d\mu}{d\mu'} =\sqrt{\frac{\det A
(t)}{\det A' (t)}} > 1,$ we find that $0<a'(t)< a(t) c (0).$
We have that:
$$
\limsup_{t \rightarrow \infty} \frac{1}{t} \int_0^t r_z'(s) \: ds  =  
\limsup_{t \rightarrow \infty} \frac{1}{t} \log a_z'(t) \leq
\lim_{t \rightarrow \infty} \frac{1}{t} \log a_z (t) =  \chi(z) $$

$$\liminf_{t \rightarrow -\infty} \frac{1}{|t|} \int_t^0 r_z'(s) \: ds =  
- \limsup_{t \rightarrow -\infty} \frac{1}{|t|} \log a_z'(t)  \geq
- \lim_{t \rightarrow -\infty} \frac{1}{|t|} \log a_z (t) =  \chi(z) ,$$

\noindent
hence 
$$\limsup_{t \rightarrow \infty} \frac{1}{t} \int_0^t r_z'(s)\: ds 
\leq \chi(z) \leq
\liminf_{t \rightarrow -\infty} \frac{1}{|t|} \int_t^0 r_z'(s) \: ds.$$

\noindent
$r_z'$ is measurable on $N,$ since continuos. Applying the following Lemma
(see \cite{BaWoEnGeo}), 
we can prove it is also integrable on $N$:

\begin{lemma}
Let $\phi^t$ be a measure preserving flow on a probability space $(X,\mu)$ and
$f: X \rightarrow \mathbb{R}$ a measurable nonnegative function; if for
almost every $x \in X$ $\limsup_{T \rightarrow + \infty} \frac{1}{T} \int_0^T 
f(\phi^t x)\:dt \leq k(x),$ where $k :X \rightarrow \mathbb{R}$ is a measurable function, 
then 
$$\int_X f(x) \:d \mu(x) \leq \int_X k(x) \:d \mu(x).$$ 
\end{lemma}

Hence, we get by Ergodic Theorem and equality of 
time averages in the future and in the past that   
$\int_N r_z(0)' d\mu = \int_N \chi(z) \, d\mu (z) = h_{\mu}(\phi).$

Finally, we use the following result (see \cite{BaWoEnGeo}):

\begin{lemma}
Given three symmetric linear operators $U,M,N$ on a 
Euclidean space such that $M$ and $N$ are nonnegative definite and $U$ is 
strictly positive definite, we get that $\traccia  [MU + NU^{-1}] \geq
2 \traccia  \sqrt{M} \sqrt{N},$ where equality holds iff $\sqrt{M}U = \sqrt{N}.$
\end{lemma}

\noindent
Since we have that $r_z'(t) = \frac{1}{2}\traccia  [V^0_{\phi^t z} - R_z^0(t){
V_{\phi^t z}^0}^{-1}],$ where $V^0_{\phi^t z}$ is (strictly) positive definite 
and $-R_z^0(t)$ is nonnegative definite, we can apply previous lemma with $U=
V_{\phi^t z}^0,$ $M=\mathbb{I}$ and $N=-R_z^0(t),$ obtaining 
$\frac{1}{2}\traccia [V^0_{\phi^t z} - R_z^0(t){V_{\phi^t z}^0}^{-1}] \geq 
\traccia  \sqrt{-R_z^0(t)},$ and hence
$$
h_{\mu}(\phi) \geq \int_N \traccia \sqrt{-R_z^0(0)} \: d\mu = \int_N \traccia 
\sqrt{-R_z(0)} \: d\mu.
$$ \hfill $\square$

\vspace{0.7cm}
\noindent
\textbf{Remark} The estimate is sharp (i.e. we have the equality) if and only if 
$V_{\phi^t z}^0 = \sqrt{-R_z^0(t)}$ for almost all $z \in N,$ which implies that 
$V^2_{\phi^t z} = -R_z(t)$ almost everywhere on $N$, and hence, by continuity, 
for every $z \in N$; this means that $\dot{V}_{\phi^t z} =0$ on $N$, i.e. 
all Jacobi curves are symmetric \cite{AgGeHam}. 

\vspace{1cm}
\noindent
\textbf{Acknowledgements}  I would like to thank Prof. A. A. Agrachev for his constant 
support and fruitful discussions.

\end{document}